# Notes on "Symmetric Bases with large 2-range"

**History**

| | | |
|---|---|---|
| 0.01 | 12-07-93 | Started |
| 0.02 | 15-07-93 | First draft completed |
| 0.03 | 26-07-93 | Results for p=22 added |
| 1.00 | 25-09-93 | Minor edits prior to sending to Bergen |
| 0.04 | 06-09-93 | Results for p=23 added |
| 1.01 | 09-07-09 | Clarification of the definitions of extensibility and symmetricisability |
| | | Add explanatory note to proof of Theorem 6 |
| | | (These additions are based on notes I made in September 2002) |
| 1.02 | 06-10-09 | Added new theorem 10 for symmetricisable p+-bases following a review of this work prompted by the publication of a paper by John Robinson on the same topic |
| 1.03 | 22-09-14 | Edited for ArxIv publication. New Abstract and Reference sections added, and all references in the body of the paper updated to reflect the new reference numbers. No other changes made. |

**Abstract**

$A_k = \{1, a_2, ... a_k\}$ is an h-basis for n if every positive integer not exceeding n can be expressed as the sum of no more than h values $a_i$; we write $n = n_h(A_k)$. An extremal h-basis $A_k$ is one for which n is as large as possible. Computing such extremal bases has become known as the Postage Stamp Problem.

This paper is inspired by and based upon a paper entitled "Symmetric bases with large 2-range for k<=75" by Svein Mossige at the University of Bergen (Mossige, Svein, [4]). Computer searches have identified some further bases which are superior to those reported in [4], and the paper also reports an improvement to one of the theoretical results.

## 1 Introduction

These notes are inspired by and based upon a paper entitled "Symmetric bases with large 2-range for k<=75" by Svein Mossige at the University of Bergen [4].

My own computer searches have identified some further bases which are superior to those reported in [4], and I have also made a small improvement (less stringent conditions) to one of the theoretical results.

In the remainder of these notes I assume that h=2; for example, I will talk of the "range $n(A_k)$" rather than the "2-range $n(2,A_k)$". In addition, the phrase "p-basis" is used to describe a special kind of h=2 basis, and the notation "6-basis" is used to mean "a p-basis for p=6" (rather than a basis for h=6).

## 2 Definitions and background

*2.1 Some notes on symmetric bases*

A basis $A_k$ is *symmetric* if $a_i + a_{k-i} = a_k$ for $1 <= i <= k-1$.

Theorem 1

If $A_k$ is symmetric and admissible, then $n(A_k) = 2a_k$.

Proof

Given $0 <= x <= a_k$, there exist $i,j$ such that $x = a_i + a_j$

$=>$ $x = (a_k - a_{k-i}) + (a_k - a_{k-j})$    - by symmetry

$=>$ $2a_k - x = a_{k-i} + a_{k-j}$

$=>$ all values $2a_k >= x >= a_k$ can also be generated.

Given an initial segment $A_j$, a symmetric set can be constructed in two ways; we use the set
$A_5 = \{1,3,4,6,11\}$ as an example:

Even $k = 2j$:

```
     1   2   1   2   5    5    2    1    2    1
  0  1   3   4   6   11   16   18   19   21   22
     └───────────────┘    └────────────────────┘
       initial segment           reflection
```

In this case, we have $a_{2j} = 2a_j$.

Odd $k = 2j-1$:

```
     1   2   1   2   5    2    1    2    1
  0  1   3   4   6   11   13   14   16   17
     └───────────────┘    └───────────────┘
       initial segment         reflection
```

In this case, we have $a_{2j-1} = a_j + a_{j-1}$.

Note that the initial segment being admissible does not guarantee that the derived symmetric basis is admissible:

Even k:

See the example above: $n(A_5) = 12$, and also $n(A_{10}) = 12$.

Odd k:

Take $A_6 = \{1,2,4,5,10,13\}$, with $n(A_6) = 15$.
Then $A_{11} = \{1,2,4,5,10,13,18,19,21,22,23\}$ also has $n(A_{11}) = 15$.

[ As it happens, the even derivation is also not admissible:
    $A_{12} = \{1,2,4,5,10,13,16,21,22,24,25,26\}$ has $n(A_{12}) = 18$. ]

## 2.2 p-bases

An admissible basis $A_{p-1}$ is called a *p-basis* if the set $\{a_i \pmod{p}: i = 1 \ldots p-1\}$ is identical to the set $\{i: i = 1 \ldots p-1\}$; in other words, the elements of the basis modulo p include each of the values 1 ... p-1 exactly once.

eg:  $A_7 = \{1, 3, 4, 5, 6, 10, 15\}$ is an 8-basis, since it is admissible ($n(A_7) = 16$), and
     $A_7 \pmod 8 = \{1, 3, 4, 5, 6, 2, 7\} = \{1, 2, 3, 4, 5, 6, 7\}$.

## 2.3 Extensible bases

A basis $A_r$ is *p-extensible* if the basis $A_{r+j}$ where:

$$A_{r+j} = \{1, a_2, \ldots, a_r=b_0, b_1, \ldots, b_j\} \text{ where } b_i = b_{i-1} + p \text{ for } 1 \le i \le j$$

is admissible for all $j \ge 0$.

Theorem 2 (below) shows that such $A_r$ must have all residues (mod p) present:

$$\{a_i \pmod{p}: i = 0 \ldots r\} = \{0, 1, 2, \ldots p-1\}$$

and so it is possible for a p-basis to be p-extensible.

## 2.4 Symmetricisable bases

Let $A_r$ be a p-extensible basis and define the symmetric basis $S(p)_k$ as follows:

$$S(p)_k = \{1, a_2, \ldots, a_r=b_0, b_1, \ldots, b_j=c_r, c_{r-1}, \ldots, c_1, c_0\} \quad \text{where } k = 2r+j$$

We say that $A_r$ is *symmetricisable* if there exists $k_0$ such that $S(p)_k$ is admissible for all $k \ge k_0$.

## 2.5 Bases with large range

Why are we interested in symmetricisable p-bases?

We're looking for bases with large range, and examination of the tables of extremal bases suggests that these are likely to be symmetric.

We know that any admissible symmetric basis $A_k$ has a range equal to $2a_k$, and so such bases are "efficient" in the sense that every generation $a_i+a_j$ contributes to the range, since all such sums are less than or equal to $2a_k$; in other words, no generation is "wasted" because its value exceeds the basis' range.

For a large range, we also require a basis that is "efficient" in the sense that as many generations as possible are unique; in other words, as few generations as possible are "wasted" because they duplicate the results of other generations.

A good starting point for bases which are efficient in this second sense is an extensible p-basis: we show later that for large enough k, the only value $b_k <= x < b_{k+1}$ which has more than one generation is the value $x = v \pmod{p}$, where $b_i = v \pmod{p}$.

So perhaps the symmetric extension of an extended p-basis may prove to be a good candidate.

In fact, such bases feature as extremal bases for k = 6 to 8, and for k = 14 to 19:

$p = 3$ $\quad A_{p-1} = \{1, 2\}$ $\quad\quad$ for k = 6 to 8
$p = 6$ $\quad A_{p-1} = \{1, 3, 4, 5, 8\}$ $\quad$ for k = 14 to 19

$\quad\quad$ eg $\quad S(6)_{14} = \{1, 3, 4, 5, 8, 14, 20, 26, 32, 35, 36, 37, 39, 40\}$

[ Perhaps surprisingly, this p=6 basis - and hence the corresponding k=14 basis - can be found in moments by hand, whereas a full computer search takes many hours! ]

We can also use such bases to construct very good bases for other values of k:

$p = 4$ $\quad A_{p-1} = \{1, 2, 3\}$

$\quad\quad S(4)_{10} = \{1, 2, 3, 7, 11, 15, 19, 20, 21, 22\}$

and $n(S(4)_{10}) = 44$, whereas the extremal range $n(10) = 46$.

The main result in [4] is a table of bases of the form $S(p)_k$ where p has been chosen to give the largest range for given k; my own results are given in section 7 below.

### 3  Extensibility

*3.1 Summary*

Let $A_{p-1}$ be a p-basis, and define $A_{p-1+j}$ as follows:

$\quad A_{p-1+j} = \{1, a_2, ..., a_{p-1}=b_0, b_1, ..., b_j\}$ $\quad$ where $b_i = b_{i-1} + p$ for $1 <= i <= j$

Mossige shows in [4] that $A_{p-1}$ is extensible if $A_{p-1+m}$ is admissible for m such that $b_{m-2} <= 2b_0 < b_{m-1}$.

We show in this section that this result can be improved in two ways:

$\quad$ i) $\quad A_{p-1}$ can be replaced by any admissible basis $A_j$ satisfying:

$\quad\quad\quad \{a_i \pmod{p}: i = 0 ... j\} = \{0, 1, 2, ..., p-1\}$ $\quad\quad$ (and hence j>=p-1)

(but not by any admissible basis that does *not* have this property; see Theorem 2).

ii) $A_j$ is extensible if and only if $A_{j+m}$ is admissible for m such that $b_m <= 2b_0 < b_{m+1}$.

I discovered these improvements by looking - unsuccessfully - for counter-examples, and the following sub-sections reflect this approach.

*3.2 Extensibility requires all residues to be present*

Theorem 2

> Let $A_{j+m} = \{1, a_2, ..., a_j=b_0, b_1, ..., b_m\}$ where $b_i = b_{i-1}+p$ for $1 <= i <= m$
>
> Then a necessary (but not sufficient) condition for $A_j$ to be extensible is that:
>
> $$\{a_i \pmod p): i = 0 ... j\} = \{0, 1, 2, ..., p-1\}$$

Proof

> "$A_j$ is extensible" means that $A_{j+m}$ is admissible for all m.
>
> Because $A_{j+m}$ is admissible, we must be able to generate all values $b_{m-1} <= x < b_m$.
>
> We choose m such that $b_{m-1} >= 2b_0$; then the generation of any value x must include at least one value greater than or equal to $b_0=a_j$:
>
> $$x = b_i + y \quad \text{for some } i >= 0$$

But we know that:

> $$x = b_{m-1} + c \quad \text{for some } 0 <= c <= p-1$$

and so we have:

> $$x = b_{m-1} + c = b_i + y$$

But $b_i = b_{m-1} \pmod p = B$, say, for all i, and so:

> $$y = c \pmod p$$

Since c runs from 0 to p-1, then so must y (mod p), which means that $A_{j+m}$ must include elements whose residues (mod p) include all possible values. $b_i \pmod p = B$ for all i, and so the remaining p-1 residues must be included in $\{0, 1, a_2, ..., a_{j-1}\}$. But $b_0 = a_j$, and so:

> $$\{a_i \pmod p): i = 0 ... j\} = \{0, 1, 2, ..., p-1\}$$

as required.

## 3.3 The first stage

Theorem 3

    Let $A_{j+m}$ be defined as in Theorem 2.

    Then $A_j$ is extensible if $A_{j+k}$ is admissible for some k such that:

$$2b_0 <= b_{k-1} \quad\quad - (1)$$

Proof

    $A_{j+k}$ is admissible => every value $b_{k-1} <= x < b_k$ can be generated, and we showed above that condition (1) means each such generation must include a value $b_i$ for some i>=0; so we have:

$$x = b_{k-1}+c = b_i+y \text{ for } 0 <= c <= p-1$$

    Adding p to both sides we have:

$$(x+p) = b_{k-1}+p+c = b_i+p+y$$

    or:

$$(x+p) = b_k+c = b_{i+1}+y \text{ for } 0 <= c <= p-1$$

    which means that all values $b_k <= (x+p) < b_{k+1}$ can be generated by $A_{j+k}$, and hence $A_{j+k+1}$ is admissible.

    Formal proof that $A_j$ is extensible follows by straightforward induction.

## 3.4 The second stage

Suppose that condition (1) in Theorem 3 above is just not met; that is, we have instead that:

$$b_{k-1} < 2b_0 <= b_k \quad\quad\quad\quad - (2)$$

The question arises as to whether bases $A_j$ exist such that $A_{j+k}$ is admissible but $A_{j+k+1}$ is not; in other words, is Theorem 3 above "sharp"?

(2) introduces the possibility that there exists a value $b_{k-1} < x < 2b_0$ whose only generations are of the form:

$$x = a_l + a_m \quad\quad \text{where } 0 <= l,m < j \quad\quad - (3)$$

If this is the case, it is easy to show that $A_{j+k+1}$ is *not* admissible:

    Suppose it is, and consider $x' = x+p > b_k >= 2b_0$; so any generation of x' must include some $b_i$ for i>0, say:

$$x' = b_i + y$$

$$\Rightarrow \quad x = x'-p = b_i-p+y = b_{i-1}+y$$

Since $i>0$, this is a generation of x which is not of the form (3) above, and so contradicts our hypothesis; therefore our assumption that $A_{j+k+1}$ is admissible is false.

In other words, the existence of an admissible basis $A_{j+k}$ with properties (2) and (3) would be sufficient to prove that Theorem 3 is "sharp"; however, a search for such a basis proved fruitless, and the following argument shows why this is the case.

Consider the value $x' = x-p$; if this has a generation of the form $b_i+y$, then x has the generation $b_{i+1}+y$ contrary to hypothesis - so x' can only have generations of the form (3). Indeed, the following picture is an inevitable consequence of our initial assumptions (2) and (3):

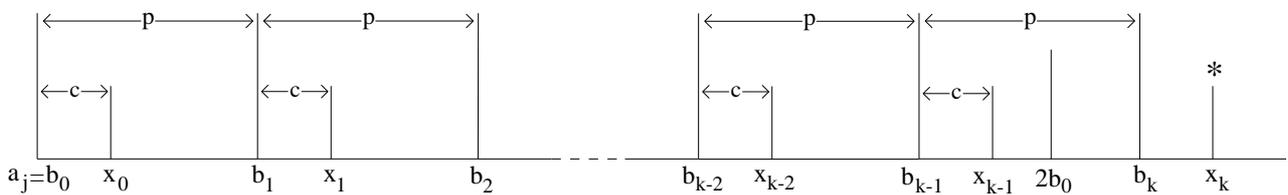

where:

a) Each $x_i$ for $0 <= i <= k-1$ has a generation involving only $\{a_i: i = 1 ... j-1\}$, but does *not* have any generation including one of $\{b_i: i = 0 ... k-1\}$.

b) $x_k$ has no generation.

We now show that condition (a) can never hold - thus proving that no basis satisfying conditions (2) and (3) can exist.

Theorem 4

Let $A_{j+m}$ be defined as in Theorem 2.

Then $A_j$ is extensible if $A_{j+k}$ is admissible for some k such that:

$$b_{k-1} < 2b_0 <= b_k$$

Proof

$A_{j+k}$ is admissible => every value $b_{k-1} <= x < b_k$ can be generated, and provided that each such generation includes a value $b_i$ for some $i>=0$ we can use the argument given in the proof of Theorem 3 to show that $A_j$ is extensible.

Suppose, however, that some value $b_{k-1} <= x < b_k$ has only generations involving $a_i$ for $i<=j-1$; then $A_{j+k}$ is a basis with the properties (2) and (3) above, and therefore must satisfy condition (a). We complete the proof by showing that condition (a) can never hold, and so no such basis can exist.

Suppose (a) holds.

We know there exists $1 <= n < j$ such that $a_n = c \pmod{p}$; consider $x = b_0+a_n$ which is the smallest generation involving some $b_i$ and $a_n$; clearly, $x = x_i \pmod{p}$.

But (a) tells us that $x$ cannot be equal to $x_i$ for any $0 <= i <= k-1$, and so we have:

$$x >= x_k$$

$$\Rightarrow \quad b_0+a_n >= b_k+c > b_k >= 2b_0$$

$$\Rightarrow \quad a_n > b_0$$

which is contrary to hypothesis.

*3.5 The final stage*

The question now naturally arises as to whether Theorem 4 is sharp. My computer searches indicated that this was not the case, but that one more refinement was possible. The result is Theorem 5, which *is* sharp; example bases which show this to be so are given in the following section.

Theorem 5

Let $A_{j+m}$ be defined as in Theorem 2.

Then $A_j$ is extensible if $A_{j+k}$ is admissible for some k such that:

$$b_k < 2b_0 <= b_{k+1}$$

Proof

Suppose the contrary: that $A_{j+k}$ is admissible, but $A_{j+k+1}$ is not.

Then there exists $b_k < x_k < b_{k+1}$ such that $x_k$ has no generation; let us write:

$$x_k = b_k+V$$

and let $a_t$ be any one of the values satisfying:

$$a_t = V \pmod{p} \quad 1 <= t <= j$$

Since $x_k$ has no generation we require that $x_k-a_t < b_0$, since otherwise a generation $x_k=b_i+a_t$ would exist for some $i>=0$. Now:

$$x_k - a_t < b_0$$

$$\Rightarrow \quad b_k + V - b_0 < a_t$$

$$\Rightarrow \quad a_t > (b_k - b_0) + V$$

$$\Rightarrow \quad a_t \geq (b_{k+1} - b_0) + V \qquad \text{- since } a_t \equiv V \pmod{p}$$

$$\Rightarrow \quad a_t \geq (2b_0 - b_0) + V = b_0 + V > b_0 = a_j \qquad \text{- which is impossible}$$

## 3.6 Examples - Theorem 5 is sharp

My computer searches were restricted to p-bases $A_{p-1}$, and I looked for bases which were at best partly extensible. In the following discussion, k is as defined in Theorem 5 above, and s is defined to be the largest value of i for which $A_{j+i}$ is admissible.

The simplest example which shows Theorem 5 to be sharp is the 5-basis $\{1, 3, 4, 7\}$:

$$A_{4+i} = \{1, 3, 4, 7, 12, 17, 22, 27, ...\}$$

We find that $A_{4+0}$ is admissible, but $A_{4+1}$ is not: $n(A_{4+0}) = 8$; so s=0. $(b_1=12) < (2b_0=14) < (b_2=17)$, and so k=1.

In fact, this is the first p-basis for which no extensions are possible.

The first p-basis for which one - but not two - extensions is possible is an 11-basis:

$$A_{10} = \{1, 2, 5, 6, 8, 9, 18, 21, 25, 26\} \text{ with extensions } \{37, 48, 59, ...\}$$

We find that $n(A_{10+1}) = 39$, so s=1; k = 2, so this is another example to show that Theorem 5 is sharp.

The first p-basis for which two - but not three - extensions is possible is a 17-basis; this also has k=s+1:

$$A_{16} = \{1, 3, 4, 7, 9, 10, 11, 14, 23, 29, 30, 32, 33, 39, 42, 53\} \text{ with extensions } \{70, 87, 104, ...\}$$

$n(A_{16+2}) = 88$, and so s=2; it is easy to see that k=3.

There are no p-bases for which three - but not infinitely many - extensions are possible for p<=21.

We can also find p-bases $A_{p-1}$ which are partly extensible, but for which s < k-1. The first such bases occur for p=10, and an example is:

$$A_9 = \{1, 3, 4, 5, 7, 8, 9, 16, 22\} \text{ with extensions } \{32, 42, 52, ...\}$$

$n(A_{9+0}) = 27$, so s=0; but k=2.

Examples where s is non-zero are harder to find; the first occurs for p=15 where we have s=1, k=3:

$A_{14}$ = {1, 2, 4, 7, 9, 12, 13, 21, 25, 29, 35, 38, 41, 48} with extensions {63, 78, 93, ... }

with $n(A_{14+1})$ = 67.

# 4 p-bases

## 4.1 Some statistics

A computer search for p-bases gives the following values of $n_p$, the number of p-bases for given p:

| p | $n_p$ | $n_p/n_{p-1}$ |
|---|---|---|
| 3 | 1 | |
| 4 | 1 | 1.00 |
| 5 | 2 | 2.00 |
| 6 | 3 | 1.50 |
| 7 | 6 | 2.00 |
| 8 | 16 | 2.67 |
| 9 | 28 | 1.75 |
| 10 | 84 | 3.00 |
| 11 | 192 | 2.29 |
| 12 | 634 | 3.30 |
| 13 | 1658 | 2.62 |
| 14 | 6277 | 3.79 |
| 15 | 18757 | 2.99 |
| 16 | 73775 | 3.93 |
| 17 | 246169 | 3.34 |
| 18 | 1044846 | 4.24 |
| 19 | 3822468 | 3.66 |
| 20 | 17365943 | 4.54 |
| 21 | 69075740 | 3.98 |
| 22 | 334698203 | 4.85 |
| 23 | 1438317540 | 4.30 |

Table 1

Note how the rate of increase fluctuates - but, overall, increases.

A quick glance at the detailed results shows that in most cases $n(A_{p-1}) < a_{p-1}+p-1$; actual counts give:

| p | n< | % | n= | % | n> | % |
|---|---|---|---|---|---|---|
| 3 | 0 | 0.0 | 1 | 100.0 | 0 | 0.0 |
| 4 | 0 | 0.0 | 1 | 100.0 | 0 | 0.0 |
| 5 | 1 | 50.0 | 1 | 50.0 | 0 | 0.0 |
| 6 | 1 | 33.3 | 2 | 66.7 | 0 | 0.0 |
| 7 | 4 | 66.7 | 2 | 33.3 | 0 | 0.0 |
| 8 | 12 | 75.0 | 3 | 18.8 | 1 | 6.2 |
| 9 | 20 | 71.4 | 6 | 21.4 | 2 | 7.1 |
| 10 | 69 | 82.1 | 11 | 13.1 | 4 | 4.8 |
| 11 | 158 | 82.3 | 23 | 12.0 | 11 | 5.7 |
| 12 | 527 | 83.1 | 54 | 8.5 | 53 | 8.4 |
| 13 | 1429 | 86.2 | 120 | 7.2 | 109 | 6.6 |
| 14 | 5495 | 87.5 | 299 | 4.8 | 483 | 7.7 |
| 15 | 16756 | 89.3 | 759 | 4.0 | 1242 | 6.6 |
| 16 | 66014 | 89.5 | 2469 | 3.3 | 5292 | 7.2 |
| 17 | 221474 | 90.0 | 7908 | 3.2 | 16787 | 6.8 |
| 18 | 941608 | 90.1 | 28764 | 2.8 | 74474 | 7.1 |

Table 2

Clearly we require $n(A_{p-1}) \geq a_{p-1}+p-1$ if $A_{p-1}$ is to be extended by even one value $b_1=a_{p-1}+p$, and so only around 10% of all p-bases are potentially extensible. We will see later that the percentage which *are* extensible is much lower, and decreases rapidly as p increases (eg 5.5% for p=16, 3.0% for p=18, 1.1% for p=21).

*4.2 Extensible p-bases eventually generate Stohr sequences*

A basis $A_j = \{1, a_2, ... , a_j\}$ defines a *Stohr sequence* $a_{j+1}, a_{j+2}, ...$ where:

$$a_{j+i+1} = n(A_{j+i}) + 1 \quad \text{for } i \geq 0$$

Reference [5] gives more details.

In general, the extension of an extensible p-basis is not necessarily its Stohr sequence; for example, consider the 8-basis {1, 3, 5, 6, 7, 10, 12}:

Stohr sequence:  1  3  5  6  7  10  12  21  23  25  34  36  38 ...
                                1  2  2  1  1  3  2  9  2  2  9  2  2 ...

Extension:       1  3  5  6  7  10  12  20  28  36  44  52  60 ...
                            1  2  2  1  1  3  2  8  8  8  8  8  8 ...

However, the following theorem shows that the Stohr sequence of a suitably extended p-basis is the same as further extension!

Theorem 6

    If $A_{p-1}$ is an extensible p-basis, then $n(A_{p-1+k}) = b_{k+1}-1$ for any k such that $b_{k+1} \geq 2b_0$.

Proof

    Since $A_{p-1}$ is extensible, we know that $n(A_{p-1+k}) \geq b_{k+1}-1$, so we have only to show that no generation exists for $b_{k+1}$ itself.

    Suppose $b_{k+1}$ has a generation; since $b_{k+1} \geq 2b_0$, such a generation must be of the form:

$$b_{k+1} = b_i + x \quad \text{for some } 0 \leq i \leq k$$

    All $b_i$ have the same residue (mod p), and so $x = 0$ (mod p); but since $A_{p-1}$ is a p-basis, the only such value x is 0*, and so there can be no generation of $b_{k+1}$.

    * Clearly none of $a_1$ to $a_{p-1} = 0$ (mod p), and so $b_0 = a_{p-1}$ - and hence all $b_i$ - also have non-zero moduli.

Corollory

    If $A_{p-1}$ is an extensible p-basis, then the Stohr sequence generated by $A_{p-1+k}$ for any k such that $b_{k+1} \geq 2b_0$ is 1-periodic with value p.

Note:

It's perfectly possible for $n(A_{p-1+k})$ to exceed $b_{k+1}-1$ for i such that $b_{k+1} < 2b_0$, as the following example for p=15 demonstrates:

$A_{14} = \{1, 2, 5, 7, 10, 11, 18, 21, 24, 27, 28, 29, 34, 38\}$ with extensions $\{53, 68, 83, ... \}$

$n(A_{14+0}) = 59 \ >= \ (b_1=53)$
$n(A_{14+1}) = 68 \ >= \ (b_2=68)$

but $b_3 > 2b_0$, and so $n(A_{14+k}) = b_{k+1}-1$ for all k>=2.

This demonstrates that Mossige's requirement (5) in [4] - that $a_{p+i} = n(A_{p+i-1})+1$ for i = 1 ... v - is not necessary; indeed, some of the bases exhibited in Table 1 in [4] do not meet this criterion (eg the optimal p=13 basis has $n(A_{12}) = a_{13}+1$).

## 5  Symmetricisability

Theorem 9

An extensible p-basis $A_{p-1}$ is symmetricisable if and only if the symmetric basis

$$S(p)_k = \{1, a_2, ... , a_{p-1}=b_0, b_1, ... , b_m=c_{p-1}, c_{p-2}, ... , c_1, c_0\} \text{ where } k = 2(p-1)+m$$

is admissible where m satisfies $b_m >= 2b_0$.

Proof

"$A_{p-1}$ is symmetricisable" means that there exists some value $m_0$ such that $S(p)_k$ is admissible for all m >= $m_0$. To prove the theorem, we show:

a)  $S(p)_k$ is admissible => $S(p)_{k+1}$ is admissible
b)  $S(p)_{k+1}$ is admissible => $S(p)_k$ is admissible

Repeated applications of (a) show that $S(p)_k$ is admissible => $A_{p-1}$ is symmetricisable, and repeated applications of (b) prove the converse.

We write:

$$S(p)_{k+1} = \{1, a_2, ... , a_{p-1}=b_0, b_1, ... , b_m, b_{m+1}=c'_{p-1}, c'_{p-2}, ... , c'_1, c'_0\}$$

and note that:

$c'_i = c_i + p$ for $0 <= i <= p-1$

a)  $S(p)_k$ is admissible => $S(p)_{k+1}$ is admissible

$A_{p-1+m+1}$ is admissible because $A_{p-1}$ is an extensible p-basis, and so all values $x' <= b_{m+1} = c'_{p-1}$ can be generated; so we only have to consider values $c'_{p-1} < x' <= c'_0$.

Now we know that the corresponding values $c_{p-1} < x \leq c_0$ have generations in $S(p)_k$ which must take one of the following forms:

$$
\begin{aligned}
x &= a_i + a_j & 0 &\leq i,j < p-1 & &- (1) \\
x &= b_i + a_j & 0 &\leq i \leq m,\ 0 \leq j < p-1 & &- (2) \\
x &= b_i + b_j & 0 &\leq i,j \leq m & &- (3) \\
x &= c_i + a_j & 0 &\leq i,j < p-1 & &- (4) \\
x &= c_i + b_j & 0 &\leq i < p-1,\ 0 \leq j \leq m & &- (5) \\
x &= c_i + c_j & 0 &\leq i,j < p-1 & &- (6)
\end{aligned}
$$

If the generation of x has one of the forms (2) to (5), we can add p to it to obtain a generation of $x' = x+p$ for $c_{p-1}+p = c'_{p-1} < x' \leq c_0+p = c'_0$:

$$
\begin{aligned}
(2) &\Rightarrow & x' &= b_i+p+a_j = b_{i+1} + a_j & 1 &\leq i+1 \leq m+1,\ 0 \leq j < p-1 \\
(3) &\Rightarrow & x' &= b_i+p+b_j = b_{i+1} + b_j & 1 &\leq i+1 \leq m+1,\ 0 \leq j \leq m \\
(4) &\Rightarrow & x' &= c_i+p+a_j = c'_i + a_j & 0 &\leq i,j < p-1 \\
(5) &\Rightarrow & x' &= c_i+p+b_j = c'_i + b_j & 0 &\leq i < p-1,\ 0 \leq j \leq m
\end{aligned}
$$

We complete this part of the proof by showing that cases (1) and (6) cannot arise:

(1): $\quad x > c_{p-1} = b_m \geq 2b_0$ and $a_i+a_j < 2a_{p-1} = 2b_0$

$\Rightarrow\quad x > a_i+a_j \quad$ for all $\ 0 \leq i,j < p-1$

(6): $\quad x \leq c_0 = c_{p-1}+a_{p-1} = b_m+b_0 \leq 2b_m$ and $c_i+c_j > 2c_{p-1} = 2b_m$

$\Rightarrow\quad x < c_i+c_j \quad$ for all $\ 0 \leq i,j < p-1$

b) $S(p)_{k+1}$ is admissible $\Rightarrow S(p)_k$ is admissible

As before, we need consider only values $c_{p-1} < x \leq c_0$, and we know that the corresponding values $c'_{p-1} < x' \leq c'_0$ have generations in $S(p)_{k+1}$ which must take one of the following forms:

$$
\begin{aligned}
x' &= a_i + a_j & 0 &\leq i,j \leq p-1 & &- (1) \\
x' &= b_i + a_j & 1 &\leq i \leq m+1,\ 0 \leq j \leq p-1 & &- (2) \\
x' &= b_i + b_j & 1 &\leq i \leq m+1,\ 1 \leq j \leq m & &- (3) \\
x' &= c'_i + a_j & 0 &\leq i,j < p-1 & &- (4) \\
x' &= c'_i + b_j & 0 &\leq i \leq p-1,\ 1 \leq j \leq m & &- (5) \\
x' &= c'_i + c'_j & 0 &\leq i,j \leq p-1 & &- (6)
\end{aligned}
$$

If the generation of x has one of the forms (2) to (5), we can subtract p from it to obtain a generation of $x = x'-p$ for $c'_{p-1}-p = c_{p-1} < x \leq c'_0-p = c_0$:

$$
\begin{aligned}
(2) &\Rightarrow & x &= b_i-p+a_j = b_{i-1} + a_j & 0 &\leq i-1 \leq m,\ 0 \leq j \leq p-1 \\
(3) &\Rightarrow & x &= b_i-p+b_j = b_{i-1} + b_j & 0 &\leq i-1 \leq m,\ 1 \leq j \leq m \\
(4) &\Rightarrow & x &= c'_i-p+a_j = c_i + a_j & 0 &\leq i,j < p-1 \\
(5) &\Rightarrow & x &= c'_i-p+b_j = c_i + b_j & 0 &\leq i \leq p-1,\ 1 \leq j \leq m
\end{aligned}
$$

On the other hand, cases (1) and (6) cannot arise:

(1):     $x' > c'_{p-1} = b_{m+1} = b_m + p >= 2b_0 + p$ and $a_i + a_j <= 2a_{p-1} = 2b_0$

=> $x' > a_i + a_j$ for all $0 <= i,j <= p-1$

(6):     $x' <= c'_0 = c'_{p-1} + a_{p-1} = b_{m+1} + b_0 < 2b_{m+1}$ and $c'_i + c'_j >= 2c'_{p-1} = 2b_{m+1}$

=> $x' < c'_i + c'_j$ for all $0 <= i,j <= p-1$

## 6 Examples and counter-examples

Let $m_0$ be defined by $b_{m_0-1} < 2b_0 <= b_{m_0}$. Then Theorem 9 characterises the behaviour of $S(p)_k$ for all $m >= m_0$ according to the behaviour at $m = m_0$: if the basis is admissible when $m=m_0$, it is admissible for all $m>m_0$; if it is not, then it is also not admissible for any $m>m_0$.

But what happens for values of $m < m_0$?

a)   Are there symmetricisable bases $A_{p-1}$ with some $S(p)_k$ for $m<m_0$ inadmissible?

b)   Are there extensible bases $A_{p-1}$ which are not symmetricisable, but for which some $S(p)_k$ for $m<m_0$ *are* admissible?

Before looking at these questions, here are some more general statistics relating to symmetricisable p-bases.

| p | $n_p$ | $n_e$ | $n_s$ | $\%_e$ | $\%_s$ |
|---|---|---|---|---|---|
| 5 | 2 | 1 | 1 | 50.0 | 100.0 |
| 6 | 3 | 2 | 2 | 66.7 | 100.0 |
| 7 | 6 | 2 | 2 | 33.3 | 100.0 |
| 8 | 16 | 4 | 4 | 25.0 | 100.0 |
| 9 | 28 | 8 | 8 | 28.6 | 100.0 |
| 10 | 84 | 15 | 15 | 17.9 | 100.0 |
| 11 | 192 | 33 | 33 | 17.2 | 100.0 |
| 12 | 634 | 99 | 99 | 15.6 | 100.0 |
| 13 | 1658 | 193 | 193 | 11.6 | 100.0 |
| 14 | 6277 | 601 | 599 | 9.6 | 99.7 |
| 15 | 18757 | 1241 | 1238 | 6.6 | 99.8 |
| 16 | 73775 | 4087 | 4062 | 5.5 | 99.4 |
| 17 | 246169 | 8883 | 8835 | 3.6 | 99.5 |
| 18 | 1044846 | 31026 | 30803 | 3.0 | 99.3 |
| 19 | 3822468 | 73367 | 72713 | 1.9 | 99.1 |
| 20 | 17365943 | 280483 | 277055 | 1.6 | 98.8 |
| 21 | 69075740 | 725490 | 715731 | 1.1 | 98.7 |

$n_p$ - number of p-bases
$n_e$ - number of extensible p-bases
$n_s$ - number of symmetricisable p-bases
$\%_e$ - $n_e/n_p$ as a percentage
$\%_s$ - $n_s/n_e$ as a percentage

Table 3

We see that the fraction of all p-bases that are extensible drops rapidly as p increases, but almost all of those that are extensible are also symmetricisable.

The first extensible p-bases which are not symmetricisable occur for p=14, and are:

$A_{13}$ = {1, 2, 4, 5, 9, 12, 13, 17, 20, 21, 22, 24, 25}
and $A_{13}$ = {1, 3, 4, 5, 8, 12, 13, 16, 20, 21, 23, 24, 25}

In each case, $m_0$=2 and $S_0$* is admissible, whereas $S_1$ and $S_2$ are not - and so the answer to (b) above is "yes" (less trivial examples follow).

[ *  In these examples, we write $S_m$ as a shorthand for $S(p)_k$ where k = 2(p-1)+m. ]

The first extensible p-bases which are not symmetricisable and for which *no* derived symmetric basis is admissible occur for p=15 and are:

$A_{14}$ = {1, 2, 5, 7, 10, 11, 18, 21, 24, 27, 28, 29, 34, 38}
and $A_{14}$ = {1, 2, 5, 7, 11, 12, 14, 18, 21, 23, 24, 25, 28, 34}

In both cases, $m_0$=3, and all $S_i$ are inadmissible for i>=0.

The first extensible p-bases which are not symmetricisable but for which some non-trivial derived symmetric basis is admissible occur for p=16 and are:

$A_{15}$ = {1, 2, 5, 8, 9, 12, 13, 19, 23, 26, 27, 30, 31, 36, 38}
and $A_{15}$ = {1, 2, 5, 8, 10, 12, 19, 22, 23, 25, 30, 31, 36, 43, 45}
and $A_{15}$ = {1, 3, 4, 8, 9, 12, 13, 18, 22, 26, 27, 30, 31, 37, 39}

In all cases, $m_0$=3.
In the first and last cases, $S_0$ and $S_1$ are admissible, but $S_2$ and $S_3$ are not.
In the second case, $S_1$ is admissible, but $S_0$, $S_2$ and $S_3$ are not.

I have also conducted a search for symmetricisable p-bases for which some $S_m$ for m<$m_0$ is *not* admissible, but so far without success: there are none for p<=21, and question (a) above remains unanswered.

# 7 Optimal bases $S_k$

## 7.1 The optimisation process

Let $S(p)_k$ be an admissible symmetric basis derived from an extended p-basis:

$$S(p)_k = \{1, a_2, \ldots, a_{p-1}=b_0, b_1, \ldots, b_j=c_{p-1}, c_{p-2}, \ldots, c_1, c_0\} \quad \text{where} \quad k = 2(p-1)+j$$

If j is even, j=2n, the basis has $k = 2(p-1)+2n$ elements, and the range is given by:

$$n(S(p)_k) = 4a_{p-1+n} = 4b_n = 4(b_0+np) = 2(2b_0+jp),$$

and if j is odd, j=2n-1, the basis has $k = 2(p-1)+2n-1$ elements, and the range is given by:

$$n(S(p)_k) = 2(a_{p-1+n}+a_{p-1+n-1}) = 2(b_n+b_{n-1}) = 2(2b_0+(2n-1)p) = 2(2b_0+jp);$$

so in both cases we have $n(S(p)_k) = 2(2b_0+jp)$ where $j = k - 2(p-1)$.

This can also be written more elegantly as $2(b_0+b_j)$, and more usefully as $2(2a_{p-1}+jp)$ - showing that for fixed p we must search for symmetricisable p-bases with largest $a_{p-1}$.

In summary, the optimisation process is as follows:

1) For each p, find the "best" symmetricisable p-bases; such p-bases will satisfy:

   a) $\{a_i \pmod{p}: i = 1 \ldots p-1\} = \{1, 2, \ldots, p-1\}$

   b) $A_{p-1}$ is admissible

   c) $A_{p-1+m} = \{1, a_2, \ldots, a_{p-1}=b_0, b_1, \ldots, b_m\}$ where $b_i = b_{i-1}+p$ for $i = 1 \ldots m$ is admissible for m such that $b_m < 2b_0 \leq b_m+p$. This guarantees that $A_{p-1}$ is extensible.

   d) $S(p)_k = \{1, a_2, \ldots, a_{p-1}=b_0, b_1, \ldots, b_j=c_{p-1}, c_{p-2}, \ldots, c_1, c_0\}$ where $k = 2(p-1)+j$ is admissible for j such that $b_{j-1} < 2b_0 \leq b_j$. This guarantees that $A_{p-1}$ is symmetricisable.

   e) $a_{p-1}$ is maximal; that is, there is no other symmetricisable p-basis with a larger $a_{p-1}$ value.

2) For given k, choose p such that $2a_{p-1}+jp$ is maximised where p and j are constrained by $k = 2(p-1)+j$.

## 7.2 Optimal solutions

Table 4 lists every symmetricisable p-basis with maximal $a_{p-1}$ for each p. We denote such a p-basis as $A^*_{p-1}$, and any symmetric basis derived from it as $S^*(p)_k$.

The table shows that the maximal value of $a_{p-1}$ does not rise smoothly with p; for example, the value for p=12 is actually less than that for p=11!

| p | $a_{p-1}$ | $A^*_{p-1}$ | | | | | | | | | | | | | | | |
|---|---|---|---|---|---|---|---|---|---|---|---|---|---|---|---|---|---|
| 5 | 4 | 1 | 2 | 3 | 4 | | | | | | | | | | | | |
| 6 | 8 | 1 | 3 | 4 | 5 | 8 | | | | | | | | | | | |
| 7 | 9 | 1 | 3 | 4 | 5 | 6 | 9 | | | | | | | | | | |
| 8 | 12 | 1 | 3 | 5 | 6 | 7 | 10 | 12 | | | | | | | | | |
| 9 | 16 | 1 | 3 | 4 | 5 | 8 | 11 | 15 | 16 | | | | | | | | |
| 10 | 19 | 1 | 2 | 5 | 6 | 8 | 13 | 14 | 17 | 19 | | | | | | | |
| 11 | 24 | 1 | 3 | 4 | 7 | 8 | 9 | 16 | 17 | 21 | 24 | | | | | | |
| 12 | 23 | 1 | 2 | 3 | 6 | 7 | 8 | 10 | 16 | 17 | 21 | 23 | | | | | |
| | | 1 | 2 | 3 | 6 | 7 | 8 | 16 | 17 | 21 | 22 | 23 | | | | | |
| | | 1 | 2 | 3 | 6 | 7 | 9 | 16 | 17 | 20 | 22 | 23 | | | | | |
| | | 1 | 2 | 3 | 7 | 8 | 10 | 16 | 17 | 18 | 21 | 23 | | | | | |
| | | 1 | 2 | 4 | 5 | 8 | 9 | 10 | 15 | 18 | 19 | 23 | | | | | |
| | | 1 | 2 | 4 | 6 | 7 | 10 | 15 | 17 | 20 | 21 | 23 | | | | | |
| | | 1 | 2 | 5 | 6 | 8 | 9 | 10 | 15 | 16 | 19 | 23 | | | | | |
| | | 1 | 2 | 5 | 6 | 8 | 10 | 15 | 16 | 19 | 21 | 23 | | | | | |
| | | 1 | 2 | 5 | 7 | 8 | 9 | 15 | 16 | 18 | 22 | 23 | | | | | |
| | | 1 | 3 | 4 | 6 | 7 | 8 | 9 | 10 | 14 | 17 | 23 | | | | | |
| | | 1 | 3 | 4 | 6 | 8 | 9 | 14 | 17 | 19 | 22 | 23 | | | | | |
| | | 1 | 3 | 4 | 6 | 9 | 10 | 14 | 17 | 19 | 20 | 23 | | | | | |
| | | 1 | 3 | 4 | 7 | 8 | 10 | 14 | 17 | 18 | 21 | 23 | | | | | |
| | | 1 | 3 | 4 | 7 | 9 | 10 | 14 | 17 | 18 | 20 | 23 | | | | | |
| 13 | 30 | 1 | 2 | 5 | 7 | 10 | 11 | 19 | 21 | 22 | 25 | 29 | 30 | | | | |
| 14 | 34 | 1 | 3 | 4 | 7 | 8 | 11 | 13 | 23 | 24 | 26 | 30 | 33 | 34 | | | |
| 15 | 41 | 1 | 2 | 5 | 6 | 8 | 9 | 13 | 19 | 22 | 27 | 29 | 33 | 40 | 41 | | |
| 16 | 45 | 1 | 2 | 3 | 7 | 8 | 9 | 12 | 15 | 22 | 26 | 30 | 36 | 37 | 43 | 45 | |
| 17 | 49 | 1 | 2 | 5 | 6 | 7 | 12 | 13 | 16 | 26 | 28 | 31 | 37 | 38 | 42 | 44 | 49 |
| 18 | 53 | 1 | 2 | 3 | 5 | 6 | 10 | 11 | 13 | 14 | 26 | 27 | 34 | 40 | 43 | 48 | 51 | 53 |
| | | 1 | 2 | 4 | 5 | 9 | 11 | 12 | 14 | 16 | 25 | 26 | 28 | 33 | 39 | 42 | 49 | 53 |
| | | 1 | 3 | 4 | 5 | 8 | 11 | 13 | 20 | 24 | 28 | 30 | 34 | 43 | 45 | 50 | 51 | 53 |
| | | 1 | 3 | 4 | 5 | 8 | 11 | 13 | 20 | 24 | 30 | 32 | 34 | 43 | 45 | 46 | 51 | 53 |
| | | 1 | 3 | 4 | 5 | 9 | 11 | 12 | 16 | 26 | 28 | 31 | 38 | 42 | 43 | 50 | 51 | 53 |
| | | 1 | 3 | 4 | 6 | 10 | 11 | 15 | 23 | 25 | 31 | 32 | 34 | 38 | 44 | 45 | 48 | 53 |
| | | 1 | 3 | 4 | 8 | 9 | 13 | 15 | 16 | 24 | 30 | 32 | 38 | 41 | 43 | 46 | 47 | 53 |
| | | 1 | 3 | 4 | 9 | 11 | 12 | 15 | 16 | 25 | 28 | 32 | 38 | 41 | 42 | 44 | 49 | 53 |
| | | 1 | 3 | 4 | 9 | 11 | 15 | 16 | 20 | 23 | 28 | 30 | 42 | 43 | 44 | 49 | 50 | 53 |
| 19 | 59 | 1 | 3 | 4 | 6 | 8 | 13 | 14 | 15 | 24 | 29 | 30 | 37 | 45 | 47 | 50 | 54 | 55 | 59 |
| | | 1 | 3 | 4 | 6 | 9 | 11 | 12 | 14 | 16 | 17 | 26 | 27 | 32 | 34 | 43 | 48 | 56 | 59 |
| | | 1 | 3 | 4 | 6 | 9 | 11 | 12 | 16 | 17 | 26 | 27 | 32 | 33 | 34 | 43 | 48 | 56 | 59 |
| | | 1 | 3 | 4 | 6 | 11 | 12 | 14 | 16 | 17 | 26 | 27 | 28 | 32 | 34 | 43 | 48 | 56 | 59 |
| | | 1 | 3 | 4 | 7 | 9 | 11 | 14 | 15 | 27 | 32 | 36 | 43 | 44 | 48 | 50 | 54 | 56 | 59 |
| | | 1 | 3 | 4 | 9 | 10 | 13 | 14 | 16 | 17 | 26 | 34 | 37 | 43 | 44 | 46 | 49 | 50 | 59 |
| | | 1 | 3 | 4 | 9 | 11 | 12 | 14 | 16 | 17 | 25 | 26 | 27 | 32 | 34 | 43 | 48 | 56 | 59 |
| | | 1 | 3 | 4 | 9 | 11 | 12 | 14 | 16 | 26 | 27 | 32 | 34 | 36 | 43 | 44 | 48 | 56 | 59 |
| | | 1 | 3 | 4 | 9 | 11 | 12 | 16 | 25 | 26 | 27 | 32 | 33 | 34 | 36 | 43 | 48 | 56 | 59 |
| | | 1 | 3 | 4 | 9 | 11 | 13 | 14 | 16 | 17 | 25 | 31 | 37 | 43 | 45 | 46 | 48 | 53 | 59 |
| | | 1 | 3 | 4 | 9 | 11 | 14 | 15 | 16 | 17 | 25 | 26 | 31 | 37 | 43 | 46 | 48 | 51 | 59 |
| 20 | 68 | 1 | 2 | 4 | 5 | 11 | 13 | 14 | 19 | 29 | 35 | 37 | 43 | 46 | 47 | 50 | 52 | 56 | 58 | 68 |
| 21 | 74 | 1 | 2 | 3 | 6 | 10 | 14 | 17 | 19 | 26 | 29 | 36 | 41 | 49 | 51 | 54 | 55 | 58 | 60 | 67 | 74 |

Table 4

Using the formula for n(S(p)$_k$) derived above, we can easily construct Table 5, from which we derive Table 6 - which, together with Table 4, corresponds to Mossige's Tables 1 and 2 in [4].

| p: | 5 | 6 | 7 | 8 | 9 | 10 | 11 | 12 | 13 | 14 | 15 | 16 | 17 | 18 | 19 | 20 | 21 |
|---|---|---|---|---|---|---|---|---|---|---|---|---|---|---|---|---|---|
| $a_{p-1}$: | 4 | 8 | 9 | 12 | 16 | 19 | 24 | 23 | 30 | 34 | 41 | 45 | 49 | 53 | 59 | 68 | 74 |
| k | | | | | | | | | | | | | | | | | |
| 8 | 16 | | | | | | | | | | | | | | | | |
| 9 | 26 | | | | | | | | | | | | | | | | |
| 10 | 36 | 32 | | | | | | | | | | | | | | | |
| 11 | 46 | 44 | | | | | | | | | | | | | | | |
| 12 | 56 | 56 | 36 | | | | | | | | | | | | | | |
| 13 | | 68 | 50 | | | | | | | | | | | | | | |
| 14 | | 80 | 64 | 48 | | | | | | | | | | | | | |
| 15 | | 92 | 78 | 64 | | | | | | | | | | | | | |
| 16 | | 104 | 92 | 80 | 64 | | | | | | | | | | | | |
| 17 | | 116 | 106 | 96 | 82 | | | | | | | | | | | | |
| 18 | | 128 | 120 | 112 | 100 | 76 | | | | | | | | | | | |
| 19 | | 140 | 134 | 128 | 118 | 96 | | | | | | | | | | | |
| 20 | | 152 | 148 | 144 | 136 | 116 | 96 | | | | | | | | | | |
| 21 | | 164 | 162 | 160 | 154 | 136 | 118 | | | | | | | | | | |
| 22 | | 176 | 176 | 176 | 172 | 156 | 140 | 92 | | | | | | | | | |
| 23 | | | | 192 | 190 | 176 | 162 | 116 | | | | | | | | | |
| 24 | | | | 208 | 208 | 196 | 184 | 140 | 120 | | | | | | | | |
| 25 | | | | | 226 | 216 | 206 | 164 | 146 | | | | | | | | |
| 26 | | | | | 244 | 236 | 228 | 188 | 172 | 136 | | | | | | | |
| 27 | | | | | 262 | 256 | 250 | 212 | 198 | 164 | | | | | | | |
| 28 | | | | | 280 | 276 | 272 | 236 | 224 | 192 | 164 | | | | | | |
| 29 | | | | | 298 | 296 | 294 | 260 | 250 | 220 | 194 | | | | | | |
| 30 | | | | | 316 | 316 | 316 | 284 | 276 | 248 | 224 | 180 | | | | | |
| 31 | | | | | | | 338 | 308 | 302 | 276 | 254 | 212 | | | | | |
| 32 | | | | | | | 360 | 332 | 328 | 304 | 284 | 244 | 196 | | | | |
| 33 | | | | | | | 382 | 356 | 354 | 332 | 314 | 276 | 230 | | | | |
| 34 | | | | | | | 404 | 380 | 380 | 360 | 344 | 308 | 264 | 212 | | | |
| 35 | | | | | | | 426 | | 406 | 388 | 374 | 340 | 298 | 248 | | | |
| 36 | | | | | | | 448 | | 432 | 416 | 404 | 372 | 332 | 284 | 236 | | |
| 37 | | | | | | | 470 | | 458 | 444 | 434 | 404 | 366 | 320 | 274 | | |
| 38 | | | | | | | 492 | | 484 | 472 | 464 | 436 | 400 | 356 | 312 | 272 | |
| 39 | | | | | | | 514 | | 510 | 500 | 494 | 468 | 434 | 392 | 350 | 312 | |
| 40 | | | | | | | 536 | | 536 | 528 | 524 | 500 | 468 | 428 | 388 | 352 | 296 |
| 41 | | | | | | | | | 562 | 556 | 554 | 532 | 502 | 464 | 426 | 392 | 338 |
| 42 | | | | | | | | | 588 | 584 | 584 | 564 | 536 | 500 | 464 | 432 | 380 |
| 43 | | | | | | | | | 614 | | 614 | 596 | 570 | 536 | 502 | 472 | 422 |
| 44 | | | | | | | | | | | 644 | 628 | 604 | 572 | 540 | 512 | 464 |
| 45 | | | | | | | | | | | 674 | 660 | 638 | 608 | 578 | 552 | 506 |
| 46 | | | | | | | | | | | 704 | 692 | 672 | 644 | 616 | 592 | 548 |
| 47 | | | | | | | | | | | 734 | 724 | 706 | 680 | 654 | 632 | 590 |
| 48 | | | | | | | | | | | 764 | 756 | 740 | 716 | 692 | 672 | 632 |
| 49 | | | | | | | | | | | 794 | 788 | 774 | 752 | 730 | 712 | 674 |
| 50 | | | | | | | | | | | 824 | 820 | 808 | 788 | 768 | 752 | 716 |
| 51 | | | | | | | | | | | 854 | 852 | 842 | 824 | 806 | 792 | 758 |
| 52 | | | | | | | | | | | 884 | 884 | 876 | 860 | 844 | 832 | 800 |
| 53 | | | | | | | | | | | | 916 | 910 | 896 | 882 | 872 | 842 |
| 54 | | | | | | | | | | | | 948 | 944 | 932 | 920 | 912 | 884 |
| 55 | | | | | | | | | | | | 980 | 978 | 968 | 958 | 952 | 926 |
| 56 | | | | | | | | | | | | 1012 | 1012 | 1004 | 996 | 992 | 968 |
| 57 | | | | | | | | | | | | | 1046 | 1040 | 1034 | 1032 | 1010 |
| 58 | | | | | | | | | | | | | 1080 | 1076 | 1072 | 1072 | 1052 |
| 59 | | | | | | | | | | | | | 1114 | 1112 | | 1112 | 1094 |
| 60 | | | | | | | | | | | | | 1148 | | | 1152 | 1136 |
| 61 | | | | | | | | | | | | | | | | 1192 | 1178 |
| 62 | | | | | | | | | | | | | | | | 1232 | 1220 |
| 63 | | | | | | | | | | | | | | | | 1272 | 1262 |
| 64 | | | | | | | | | | | | | | | | 1312 | 1304 |
| 65 | | | | | | | | | | | | | | | | 1352 | 1346 |
| 66 | | | | | | | | | | | | | | | | 1392 | 1388 |
| 67 | | | | | | | | | | | | | | | | 1432 | 1430 |
| 68 | | | | | | | | | | | | | | | | 1472 | 1472 |
| 69 | | | | | | | | | | | | | | | | | 1514 |
| 70 | | | | | | | | | | | | | | | | | 1556 |
| 71 | | | | | | | | | | | | | | | | | 1598 |
| 72 | | | | | | | | | | | | | | | | | 1640 |
| 73 | | | | | | | | | | | | | | | | | 1682 |
| 74 | | | | | | | | | | | | | | | | | 1724 |
| 75 | | | | | | | | | | | | | | | | | 1766 |

This table gives $n(S^*(p)_k)$ for different values of p and k

Table 5

| $k_{min}$ | $k_{max}$ | range | p |
|---|---|---|---|
| 8 | 12 | 16 | 5 |
| 12 | 22 | 56 | 6 |
| 22 | 22 | 176 | 7 |
| 22 | 24 | 176 | 8 |
| 24 | 30 | 208 | 9 |
| 30 | 30 | 316 | 10 |
| 30 | 40 | 316 | 11 |
| 40 | 43 | 536 | 13 |
| 43 | 52 | 614 | 15 |
| 52 | 56 | 884 | 16 |
| 56 | 59 | 1012 | 17 |
| 60 | 68 | 1152 | 20 |
| 68 | | 1472 | 21 |

Table 6

As an example, the second line in Table 6 means that the optimal bases $S^*(p)_k$ for k = 12 to 22 inclusive are derived from the maximal 6-basis $A^*_5$, and the range for the minimum value of k (=12) is 56; the range for other values is equal to 56 + 2(k-12)p = 56 + 12(k-12).

    eg    The basis for k=14 has range = 56+24 = 80, and is:

$$S^*(6)_{14} = \{1, \ 3, \ 4, \ 5, \ 8, 14, 20, 26, 32, 35, 36, 37, 39, 40\}$$
$$\phantom{S^*(6)_{14} = \{}1 \ \ 2 \ \ 1 \ \ 1 \ \ 3 \ \ 6 \ \ 6 \ \ 6 \ \ 6 \ \ 3 \ \ 1 \ \ 1 \ \ 2 \ \ 1$$

Table 6 differs from Table 1 in [4] as follows:

- Mossige's p=8 basis with $a_p$ different from $a_{p-1}+p$ is omitted (see below)
- I have included improved bases for p=16 and 21, and results for p=20 are included

and although p=19 no longer features in the "best" list, Table 4 above includes two new $A^*_{18}$ which are not present in Tables 1 and 2 of [4].

*7.3 What happens if $a_p$ is allowed to vary*

Mossige [4] observed that for p=8 an improved basis appears in which $a_p$ differs from $a_{p-1}+p$. To investigate this further, I took every p-basis, extended it first by each possible $a_p$ satisfying:

$$a_{p-1} < a_p <= n(A_{p-1})+1$$

and then determined whether the resulting basis was both extensible and symmetricisable. I shall denote such a basis a *p+-basis*, and the maximal p+-bases found are listed in Table 7; they exhibit the following improvements over Table 2:

- The calculations have been extended to include p=22 and p=23,
- The p+-bases for p=8 and p=19 have improved optimal values "$a_{p-1}$"*.
- There are extra p+-bases for p=5, 10, 12, 14 and 18 which are as good as the corresponding optimal p-bases.

[ * "$a_{p-1}$" is defined as $a_p$-p (regardless of the true value of $a_{p-1}$) for comparison with Table 4. ]

The following theorem is analogous to Theorem 9 (which dealt with symmetricisable p-bases).

Theorem 10

The extensible p+-basis $A_p = \{1, a_2, \ldots a_{p-1}, a_p\}$ is symmetricisable if the symmetric basis:

$$S(p)_k = \{1, a_2, \ldots, a_p = b_0, b_1, \ldots, b_m = c_p, c_{p-1}, \ldots, c_1, c_0\}, \quad k = 2p + m$$

where $b_i = b_{i-1} + p$ for $1 \leq i \leq m$

is admissible for some m such that $b_m \geq 2b_0$.

Proof

Consider:

$$S(p)_{k+1} = \{1, a_2, \ldots, a_p = b_0, b_1, \ldots, b_{m+1} = c'_p, c'_{p-1}, \ldots, c'_1, c'_0\}$$

We shall show that $S(p)_k$ is admissible $\Rightarrow S(p)_{k+1}$ is admissible, and hence, by induction, that $A_p$ is symmetricisable.

Since $A_p$ is extensible, we know that we can generate all values $\leq c'_p$, and so we have only to consider values $c'_p < x' < c'_0$. It is clear that $c'_i = c_i + p$, and so because $S(p)_k$ is admissible we know that $x = x' - p$, $c_p < x < c_0$ has a generation. If such a generation includes at least one of $b_i$ or $c_i$ then we can derive a generation for x' by replacing $b_i$ by $b_{i+1}$, or $c_i$ by $c'_i$:

$$x = b_i + y \implies x' = b_{i+1} + y \quad \text{for } 0 \leq i \leq m$$
$$x = c_i + y \implies x' = c'_i + y$$

But if the generation of x has neither $b_i$ nor $c_i$ then it must be of the form:

$$x = a_i + a_j \quad \text{for } i, j \leq p$$

and so $x \leq 2a_p = 2b_0 \leq b_m = c_p$ - contrary to our hypothesis that $c_p < x < c_0$. So each generation of x must include either $b_i$ or $c_i$.

| p | $a_{p-1}$ | $A^*_{p-1}$ |
|---|---|---|
| 5 | 4 | 1 2 3 4 9 |
|   |   | 1 3 4 7 9 |
| 6 | 8 | 1 3 4 5 8 14 |
| 7 | 9 | 1 3 4 5 6 9 16 |
| 8 | 13 | 1 3 4 6 10 13 15 21 |
| 9 | 16 | 1 3 4 5 8 11 15 16 25 |
| 10 | 19 | 1 2 5 6 8 13 14 17 19 29 |
|   |   | 1 3 4 6 8 12 17 19 25 29 |
| 11 | 24 | 1 3 4 7 8 9 16 17 21 24 35 |
| 12 | 23 | 1 2 3 5 7 8 16 18 22 23 33 35 |
|   |   | 1 2 3 6 7 8 10 16 17 21 23 35 |
|   |   | 1 2 3 6 7 8 16 17 21 22 23 35 |
|   |   | 1 2 3 6 7 9 10 20 23 28 29 35 |
|   |   | 1 2 3 6 7 9 16 17 20 22 23 35 |
|   |   | 1 2 3 7 8 10 16 17 18 21 23 35 |
|   |   | 1 2 4 5 8 9 10 15 18 19 23 35 |
|   |   | 1 2 4 5 8 10 15 21 23 30 31 35 |
|   |   | 1 2 4 6 7 9 10 20 23 27 29 35 |
|   |   | 1 2 4 6 7 10 15 17 20 21 23 35 |
|   |   | 1 2 4 6 8 9 17 19 22 23 27 35 |
|   |   | 1 2 4 6 9 10 15 19 20 23 29 35 |
|   |   | 1 2 4 6 9 10 17 19 20 23 27 35 |
|   |   | 1 2 5 6 8 9 10 15 16 19 23 35 |
|   |   | 1 2 5 6 8 10 15 16 19 21 23 35 |
|   |   | 1 2 5 7 8 9 15 16 18 22 23 35 |
|   |   | 1 2 5 7 8 10 18 21 23 27 28 35 |
|   |   | 1 3 4 5 6 10 14 20 21 23 31 35 |
|   |   | 1 3 4 5 8 10 14 21 23 30 31 35 |
|   |   | 1 3 4 6 7 8 9 10 14 17 23 35 |
|   |   | 1 3 4 6 8 9 14 17 19 22 23 35 |
|   |   | 1 3 4 6 9 10 14 17 19 20 23 35 |
|   |   | 1 3 4 6 10 14 19 20 21 23 29 35 |
|   |   | 1 3 4 7 8 10 14 17 18 21 23 35 |
|   |   | 1 3 4 7 8 10 18 23 26 29 33 35 |
|   |   | 1 3 4 7 9 10 14 17 18 20 23 35 |
|   |   | 1 3 4 7 9 14 18 20 22 23 29 35 |
|   |   | 1 3 5 6 9 10 16 20 23 26 31 35 |
|   |   | 1 3 5 7 8 10 16 18 21 23 26 35 |
| 13 | 30 | 1 2 5 7 10 11 19 21 22 25 29 30 43 |
| 14 | 34 | 1 3 4 7 8 11 13 23 24 26 30 33 34 48 |
|   |   | 1 3 5 7 8 10 11 18 23 27 30 34 40 48 |
| 15 | 41 | 1 2 5 6 8 9 13 19 22 27 29 33 40 41 56 |
| 16 | 45 | 1 2 3 7 8 9 12 15 22 26 30 36 37 43 45 61 |
| 17 | 49 | 1 2 5 6 7 12 13 16 26 28 31 37 38 42 44 49 66 |
| 18 | 53 | 1 2 3 5 6 10 11 13 14 26 27 34 40 43 48 51 53 71 |
|   |   | 1 2 3 7 10 12 14 16 23 24 29 33 40 44 49 53 63 71 |
|   |   | 1 2 4 5 9 11 12 14 16 25 26 28 33 39 42 49 53 71 |
|   |   | 1 2 4 5 9 11 12 15 24 26 31 32 34 43 46 53 57 71 |
|   |   | 1 2 5 7 11 14 15 16 22 24 27 30 44 46 49 53 57 71 |
|   |   | 1 2 5 8 9 12 14 15 24 25 29 31 34 40 46 53 57 71 |
|   |   | 1 3 4 5 8 11 13 20 24 28 30 34 43 45 50 51 53 71 |
|   |   | 1 3 4 5 8 11 13 20 24 30 32 34 43 45 46 51 53 71 |
|   |   | 1 3 4 5 9 11 12 16 26 28 31 38 42 43 50 51 53 71 |
|   |   | 1 3 4 6 10 11 15 23 25 31 32 34 38 44 45 48 53 71 |
|   |   | 1 3 4 8 9 13 15 16 24 30 32 38 41 43 46 47 53 71 |
|   |   | 1 3 4 9 11 12 15 16 25 28 32 38 41 42 44 49 53 71 |
|   |   | 1 3 4 9 11 15 16 20 23 28 30 42 43 44 49 50 53 71 |
| 19 | 61 | 1 2 3 6 9 11 12 15 16 27 32 37 45 48 52 55 61 62 80 |
| 20 | 68 | 1 2 4 5 11 13 14 19 29 35 37 43 46 47 50 52 56 58 68 88 |
| 21 | 74 | 1 2 3 6 10 14 17 19 26 29 36 41 49 51 54 55 58 60 67 74 95 |
| 22 | 83 | 1 3 5 7 8 12 14 18 26 32 33 42 43 50 60 63 68 79 81 83 97 105 |
| 23 | 86 | 1 2 5 6 8 10 14 15 26 32 34 43 50 53 59 62 64 65 81 86 90 91 109 |

Note that the column headed $a_{p-1}$ contains the value $a_p - p$ for comparison with Table 4 above.

Table 7

Tables 8 and 9 correspond to Tables 5 and 6 respectively.

| p: | 5 | 6 | 7 | 8 | 9 | 10 | 11 | 12 | 13 | 14 | 15 | 16 | 17 | 18 | 19 | 20 | 21 | 22 | 23 |
|---|---|---|---|---|---|---|---|---|---|---|---|---|---|---|---|---|---|---|---|
| $a_{p-1}$: | 4 | 8 | 9 | 13 | 16 | 19 | 24 | 23 | 30 | 34 | 41 | 45 | 49 | 53 | 61 | 68 | 74 | 83 | 86 |
| k | | | | | | | | | | | | | | | | | | | |
| 8 | 16 | | | | | | | | | | | | | | | | | | |
| 9 | 26 | | | | | | | | | | | | | | | | | | |
| 10 | 36 | 32 | | | | | | | | | | | | | | | | | |
| 11 | 46 | 44 | | | | | | | | | | | | | | | | | |
| 12 | 56 | 56 | 36 | | | | | | | | | | | | | | | | |
| 13 | | 68 | 50 | | | | | | | | | | | | | | | | |
| 14 | | 80 | 64 | 52 | | | | | | | | | | | | | | | |
| 15 | | 92 | 78 | 68 | | | | | | | | | | | | | | | |
| 16 | | 104 | 92 | 84 | 64 | | | | | | | | | | | | | | |
| 17 | | 116 | 106 | 100 | 82 | | | | | | | | | | | | | | |
| 18 | | 128 | 120 | 116 | 100 | 76 | | | | | | | | | | | | | |
| 19 | | 140 | 134 | 132 | 118 | 96 | | | | | | | | | | | | | |
| 20 | | 152 | 148 | 148 | 136 | 116 | 96 | | | | | | | | | | | | |
| 21 | | 164 | | 164 | 154 | 136 | 118 | | | | | | | | | | | | |
| 22 | | | | 180 | 172 | 156 | 140 | 92 | | | | | | | | | | | |
| 23 | | | | 196 | 190 | 176 | 162 | 116 | | | | | | | | | | | |
| 24 | | | | 212 | 208 | 196 | 184 | 140 | 120 | | | | | | | | | | |
| 25 | | | | 228 | 226 | 216 | 206 | 164 | 146 | | | | | | | | | | |
| 26 | | | | 244 | 244 | 236 | 228 | 188 | 172 | 136 | | | | | | | | | |
| 27 | | | | | 262 | 256 | 250 | 212 | 198 | 164 | | | | | | | | | |
| 28 | | | | | 280 | 276 | 272 | 236 | 224 | 192 | 164 | | | | | | | | |
| 29 | | | | | 298 | 296 | 294 | 260 | 250 | 220 | 194 | | | | | | | | |
| 30 | | | | | 316 | 316 | 316 | 284 | 276 | 248 | 224 | 180 | | | | | | | |
| 31 | | | | | | | 338 | 308 | 302 | 276 | 254 | 212 | | | | | | | |
| 32 | | | | | | | 360 | 332 | 328 | 304 | 284 | 244 | 196 | | | | | | |
| 33 | | | | | | | 382 | 356 | 354 | 332 | 314 | 276 | 230 | | | | | | |
| 34 | | | | | | | 404 | 380 | 380 | 360 | 344 | 308 | 264 | 212 | | | | | |
| 35 | | | | | | | 426 | | 406 | 388 | 374 | 340 | 298 | 248 | | | | | |
| 36 | | | | | | | 448 | | 432 | 416 | 404 | 372 | 332 | 284 | 244 | | | | |
| 37 | | | | | | | 470 | | 458 | 444 | 434 | 404 | 366 | 320 | 282 | | | | |
| 38 | | | | | | | 492 | | 484 | 472 | 464 | 436 | 400 | 356 | 320 | 272 | | | |
| 39 | | | | | | | 514 | | 510 | 500 | 494 | 468 | 434 | 392 | 358 | 312 | | | |
| 40 | | | | | | | 536 | | 536 | 528 | 524 | 500 | 468 | 428 | 396 | 352 | 296 | | |
| 41 | | | | | | | | | 562 | 556 | 554 | 532 | 502 | 464 | 434 | 392 | 338 | | |
| 42 | | | | | | | | | 588 | 584 | 584 | 564 | 536 | 500 | 472 | 432 | 380 | 332 | |
| 43 | | | | | | | | | 614 | | 614 | 596 | 570 | 536 | 510 | 472 | 422 | 376 | |
| 44 | | | | | | | | | | | 644 | 628 | 604 | 572 | 548 | 512 | 464 | 420 | 344 |
| 45 | | | | | | | | | | | 674 | 660 | 638 | 608 | 586 | 552 | 506 | 464 | 390 |
| 46 | | | | | | | | | | | 704 | 692 | 672 | 644 | 624 | 592 | 548 | 508 | 436 |
| 47 | | | | | | | | | | | 734 | 724 | 706 | 680 | 662 | 632 | 590 | 552 | 482 |
| 48 | | | | | | | | | | | 764 | 756 | 740 | 716 | 700 | 672 | 632 | 596 | 528 |
| 49 | | | | | | | | | | | 794 | 788 | 774 | 752 | 738 | 712 | 674 | 640 | 574 |
| 50 | | | | | | | | | | | 824 | 820 | 808 | 788 | 776 | 752 | 716 | 684 | 620 |
| 51 | | | | | | | | | | | 854 | 852 | 842 | 824 | 814 | 792 | 758 | 728 | 666 |
| 52 | | | | | | | | | | | 884 | 884 | 876 | 860 | 852 | 832 | 800 | 772 | 712 |
| 53 | | | | | | | | | | | | 916 | 910 | 896 | 890 | 872 | 842 | 816 | 758 |
| 54 | | | | | | | | | | | | 948 | 944 | 932 | 928 | 912 | 884 | 860 | 804 |
| 55 | | | | | | | | | | | | 980 | 978 | 968 | 966 | 952 | 926 | 904 | 850 |
| 56 | | | | | | | | | | | | 1012 | 1012 | 1004 | 1004 | 992 | 968 | 948 | 896 |
| 57 | | | | | | | | | | | | | 1046 | | 1042 | 1032 | 1010 | 992 | 942 |
| 58 | | | | | | | | | | | | | 1080 | | 1080 | 1072 | 1052 | 1036 | 988 |
| 59 | | | | | | | | | | | | | | | 1118 | 1112 | 1094 | 1080 | 1034 |
| 60 | | | | | | | | | | | | | | | 1156 | 1152 | 1136 | 1124 | 1080 |
| 61 | | | | | | | | | | | | | | | 1194 | 1192 | 1178 | 1168 | 1126 |
| 62 | | | | | | | | | | | | | | | 1232 | 1232 | 1220 | 1212 | 1172 |
| 63 | | | | | | | | | | | | | | | | 1272 | 1262 | 1256 | 1218 |
| 64 | | | | | | | | | | | | | | | | 1312 | 1304 | 1300 | 1264 |
| 65 | | | | | | | | | | | | | | | | 1352 | 1346 | 1344 | 1310 |
| 66 | | | | | | | | | | | | | | | | 1392 | 1388 | 1388 | 1356 |
| 67 | | | | | | | | | | | | | | | | 1432 | | 1432 | 1402 |
| 68 | | | | | | | | | | | | | | | | | | 1476 | 1448 |
| 69 | | | | | | | | | | | | | | | | | | 1520 | 1494 |
| 70 | | | | | | | | | | | | | | | | | | 1564 | 1540 |
| 71 | | | | | | | | | | | | | | | | | | 1608 | 1586 |
| 72 | | | | | | | | | | | | | | | | | | 1652 | 1632 |
| 73 | | | | | | | | | | | | | | | | | | 1696 | 1678 |
| 74 | | | | | | | | | | | | | | | | | | 1740 | 1724 |
| 75 | | | | | | | | | | | | | | | | | | 1784 | 1770 |
| 76 | | | | | | | | | | | | | | | | | | 1828 | 1816 |
| 77 | | | | | | | | | | | | | | | | | | 1872 | 1862 |
| 78 | | | | | | | | | | | | | | | | | | 1916 | 1908 |
| 79 | | | | | | | | | | | | | | | | | | 1960 | 1954 |
| 80 | | | | | | | | | | | | | | | | | | 2004 | 2000 |
| 81 | | | | | | | | | | | | | | | | | | 2048 | 2046 |
| 82 | | | | | | | | | | | | | | | | | | 2092 | 2092 |
| 83 | | | | | | | | | | | | | | | | | | | 2138 |

Table 8

| $k_{min}$ | $k_{max}$ | range | p |
|---|---|---|---|
| 8 | 12 | 16 | 5 |
| 12 | 21 | 56 | 6 |
| 21 | 26 | 164 | 8 |
| 26 | 30 | 244 | 9 |
| 30 | 30 | 316 | 10 |
| 30 | 40 | 316 | 11 |
| 40 | 43 | 536 | 13 |
| 43 | 52 | 614 | 15 |
| 52 | 56 | 884 | 16 |
| 56 | 58 | 1012 | 17 |
| 58 | 62 | 1080 | 19 |
| 62 | 67 | 1232 | 20 |
| 67 | 82 | 1432 | 22 |
| 82 | | 2092 | 23 |

Table 9

*7.4 What if $a_{p+1}$, $a_{p+2}$, ... are allowed to vary?*

In a p+-basis, we allow any admissible value for $a_p$ before considering extension, and it is natural to ask whether a "p++-basis" - in which any admissible values for both $a_p$ and $a_{p+1}$ are allowed - might do even better; however my computer search shows there are no such improved bases (nor any new p++-bases with the same maximal "$a_{p-1}$") for any p<=21.

If we allow sufficient '+'s, we can always turn *any* p-basis into an extensible p++..+-basis*; it is interesting to speculate whether bases derived in this way will eventually appear as the generators of optimal $S^*(p)_k$ for large k.

  [ * Let A = {1, $a_2$, ... , $a_{p-1}$=$b_0$, $b_1$, ... , $b_m$}  where  $b_{m-1}$< $2b_0$ <= $b_m$  and  $b_i$ = $b_{i-1}$ + p  for i>0.

   We know that if A is admissible, then A is extensible.

   Suppose that A is *not* admissible, with values C = {$c_1$, $c_2$, ... , $c_n$},  $b_0$ < $c_i$ < $b_m$, which cannot be generated.

   Then define A' as the union of A and C; clearly A' is admissible, and by Theorem 4 above it is extensible. ]

How could such a basis arise?

Suppose that for a given value of p the optimal symmetricisable p-basis has $a_{p-1}$ = X; this means that no p-basis with $a_{p-1}$>X is symmetricisable. On the other hand, there might exist a p+-basis with $a_{p-1}$>X which is symmetricisable; if it also turns out that $a_p$ > X+p then this p+-basis is superior to all p-bases.

This is what happens with p=8 (cf Tables 4 and 7), where X=12:

The best symmetricisable p-basis extends as follows:

$$1, 3, 5, 6, 7, 10, 12, 20, 28, 36, 44 \ldots$$

But there is a better p+-basis with $a_{p-1}>X$ that extends as follows:

$$1, 3, 4, 6, 10, 13, 15, 21, 29, 37, 45 \ldots$$

In fact, in all cases where p+-bases in Table 7 exceed or equal the p-bases of Table 4, we find that $a_{p-1}>X$; we might expect this to be the case, since this allows $a_p<a_{p-1}+p$ while still allowing the possibility that $a_p>X+p$.

On the other hand, a p-basis with $a_{p-1}<X$ could also turn into a good p+-basis if its range is greater than $a_{p-1}+p$; so there may be scope here for improved bases for large k, although none has yet been found.

In either case, the general requirement is that:

$$(a_q - a_{p-1}) > X + (q-p+1)p$$

for a symmetricisable p++...+-basis $\{1, a_2, \ldots, a_{p-1}, a_p, \ldots, a_q\}$, and there is no obvious reason why such bases will not exist for large p.

## 8  Some statistics

Table 10 is an analysis of all p-bases for p=21 according to the value of $a_{p-1}$, and is typical of similar analyses for other values of p.

| $a_{p-1}$ | $n_p$ | $n_e$ | $n_s$ | $a_{p-1}$ | $n_p$ | $n_e$ | $n_s$ | $a_{p-1}$ | $n_p$ | $n_e$ | $n_s$ |
|---|---|---|---|---|---|---|---|---|---|---|---|
| 111 | 1 | 0 | 0 | 80 | 1291981 | 0 | 0 | 49 | 240095 | 50790 | 50547 |
| 110 | 1 | 0 | 0 | 79 | 1596027 | 0 | 0 | 48 | 190903 | 50064 | 49909 |
| 109 | 0 | 0 | 0 | 78 | 1783697 | 0 | 0 | 47 | 114720 | 35101 | 35045 |
| 108 | 1 | 0 | 0 | 77 | 2078737 | 0 | 0 | 46 | 100561 | 35914 | 35867 |
| 107 | 4 | 0 | 0 | 76 | 2256038 | 0 | 0 | 45 | 56187 | 22044 | 22020 |
| 106 | 0 | 0 | 0 | 75 | 2561807 | 0 | 0 | 44 | 57364 | 25596 | 25571 |
| 105 | 0 | 0 | 0 | 74 | 2613923 | 1 | 1 | 43 | 0 | 0 | 0 |
| 104 | 31 | 0 | 0 | 73 | 2857064 | 0 | 0 | 42 | 0 | 0 | 0 |
| 103 | 98 | 0 | 0 | 72 | 2820304 | 0 | 0 | 41 | 102618 | 55663 | 55640 |
| 102 | 148 | 0 | 0 | 71 | 2987439 | 2 | 2 | 40 | 55562 | 35901 | 35883 |
| 101 | 291 | 0 | 0 | 70 | 2798924 | 7 | 7 | 39 | 29629 | 21783 | 21781 |
| 100 | 462 | 0 | 0 | 69 | 2905337 | 16 | 16 | 38 | 15274 | 12493 | 12489 |
| 99 | 987 | 0 | 0 | 68 | 2289365 | 43 | 41 | 37 | 7891 | 6907 | 6906 |
| 98 | 1467 | 0 | 0 | 67 | 2393590 | 49 | 40 | 36 | 3970 | 3632 | 3632 |
| 97 | 2875 | 0 | 0 | 66 | 1738165 | 123 | 99 | 35 | 2040 | 1925 | 1925 |
| 96 | 4269 | 0 | 0 | 65 | 2027028 | 476 | 440 | 34 | 1016 | 977 | 977 |
| 95 | 7264 | 0 | 0 | 64 | 0 | 0 | 0 | 33 | 526 | 513 | 513 |
| 94 | 11211 | 0 | 0 | 63 | 0 | 0 | 0 | 32 | 258 | 254 | 254 |
| 93 | 17872 | 0 | 0 | 62 | 5240097 | 3745 | 3357 | 31 | 136 | 135 | 135 |
| 92 | 25620 | 0 | 0 | 61 | 4610884 | 5213 | 4745 | 30 | 66 | 66 | 66 |
| 91 | 37958 | 0 | 0 | 60 | 3830769 | 8979 | 8337 | 29 | 36 | 36 | 36 |
| 90 | 51261 | 0 | 0 | 59 | 3222152 | 11684 | 11029 | 28 | 17 | 17 | 17 |
| 89 | 58609 | 0 | 0 | 58 | 2623745 | 17577 | 16642 | 27 | 10 | 10 | 10 |
| 88 | 76979 | 0 | 0 | 57 | 2128773 | 22997 | 22008 | 26 | 4 | 4 | 4 |
| 87 | 80716 | 0 | 0 | 56 | 1687881 | 28896 | 27930 | 25 | 3 | 3 | 3 |
| 86 | 134481 | 0 | 0 | 55 | 1323339 | 31599 | 30685 | 24 | 1 | 1 | 1 |
| 85 | 0 | 0 | 0 | 54 | 1039183 | 38812 | 37965 | 23 | 1 | 1 | 1 |
| 84 | 0 | 0 | 0 | 53 | 780788 | 44482 | 43632 | 22 | 0 | 0 | 0 |
| 83 | 714741 | 0 | 0 | 52 | 604270 | 44420 | 43804 | 21 | 0 | 0 | 0 |
| 82 | 861192 | 0 | 0 | 51 | 443672 | 50399 | 49950 | 20 | 1 | 1 | 1 |
| 81 | 1131295 | 0 | 0 | 50 | 342038 | 56139 | 55768 | | | | |

Table 10

We immediately see that:

- Most extensible p-bases are symmetricisable.

- Even the best extensible p-bases are usually symmetricisable (indeed, for all p<=21, at least one of the extensible p-bases with maximum $a_{p-1}$ is also symmetricisable).

- The best p-bases are not extensible.

It is this last point which suggests that p++..+-bases may exist which are superior to p-bases and p+-bases for larger k: there are many potential candidates with $a_{p-1}$ values much greater than X.

Table 11 examines this point more closely, showing for each p:

    $v_1$   - the maximum $a_{p-1}$ of any p-basis
    $v_2$   - the maximum $a_{p-1}$ of any extensible p-basis (ie "X")

We see that while the ratio $v_1/p$ increases as p increases, the ratio $v_2/v_1$ remians remarkably constant, being very close to two thirds.

| p | $v_1$ | $v_2$ | $v_1/p$ | $v_2/v_1$ |
|---|---|---|---|---|
| 5 | 7 | 4 | 1.40 | 0.57 |
| 6 | 9 | 8 | 1.50 | 0.89 |
| 7 | 13 | 9 | 1.86 | 0.69 |
| 8 | 15 | 12 | 1.87 | 0.80 |
| 9 | 21 | 16 | 2.33 | 0.76 |
| 10 | 26 | 19 | 2.60 | 0.73 |
| 11 | 30 | 24 | 2.73 | 0.80 |
| 12 | 35 | 23 | 2.92 | 0.66 |
| 13 | 44 | 30 | 3.38 | 0.68 |
| 14 | 51 | 34 | 3.64 | 0.67 |
| 15 | 58 | 41 | 3.87 | 0.71 |
| 16 | 63 | 45 | 3.94 | 0.71 |
| 17 | 73 | 49 | 4.29 | 0.67 |
| 18 | 83 | 53 | 4.61 | 0.64 |
| 19 | 91 | 59 | 4.79 | 0.65 |
| 20 | 99 | 68 | 4.95 | 0.69 |
| 21 | 111 | 74 | 5.29 | 0.67 |

Table 11

The remaining diagrams show various distributions for p=21 as bar charts.

Table 12 shows the number of p-bases plotted against $a_{p-1}$; the gaps are easily explained:

- The "smallest" p-basis is $\{1, 2, \ldots, p-1\}$, and so there is no p-basis with $a_{p-1} < p-1$.

- Every p-basis must include $a_0=0$ and $a_1=1$, so no values of $a_{p-1}$ are ever equivalent to 0 (mod p) or 1 (mod p).

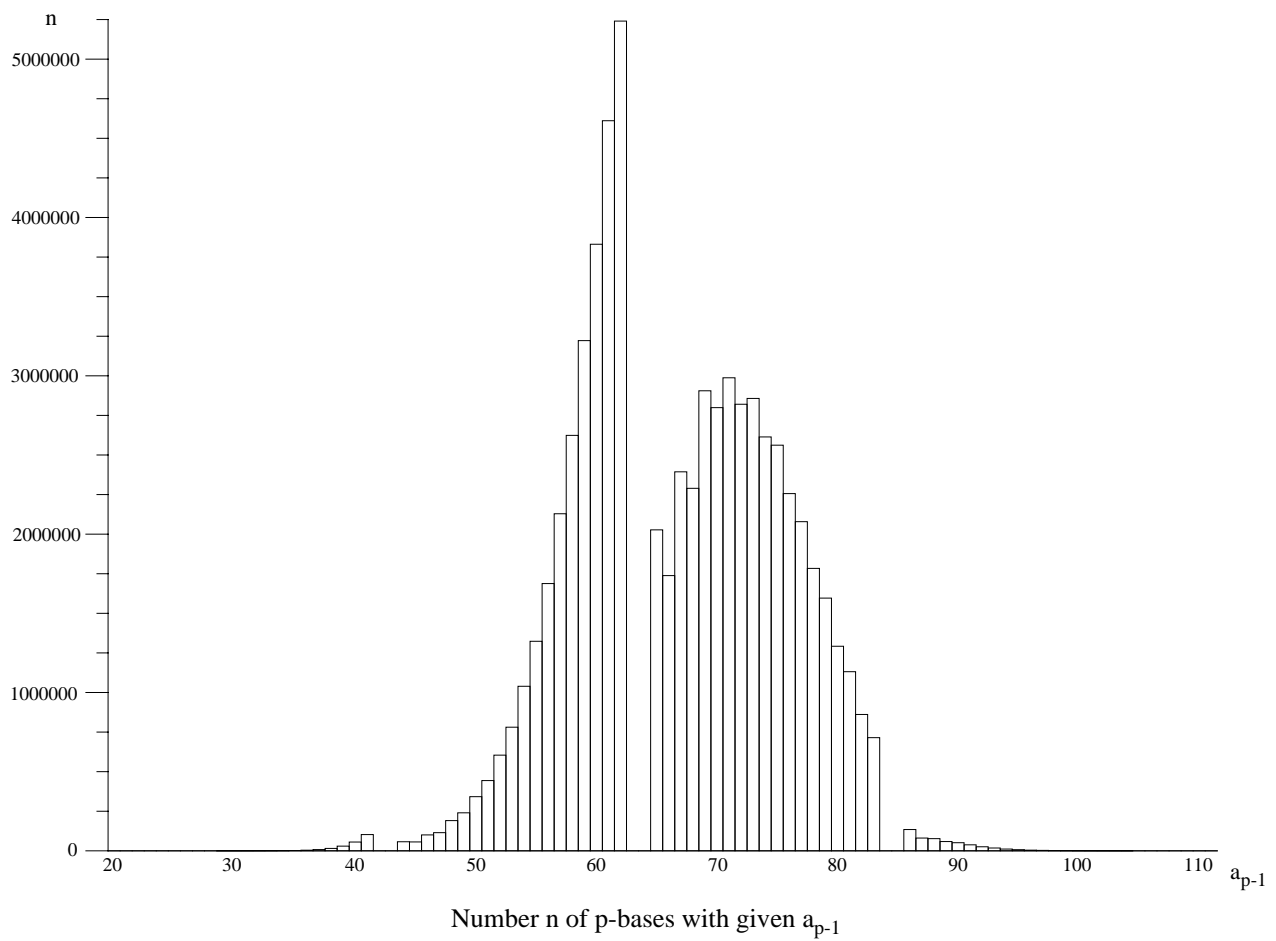

Number n of p-bases with given $a_{p-1}$

Table 12

Table 13 shows the number of extensible p-bases plotted against $a_{p-1}$; the number of these which are not also symmetricisable is indicated by the small section at the top of each bar.

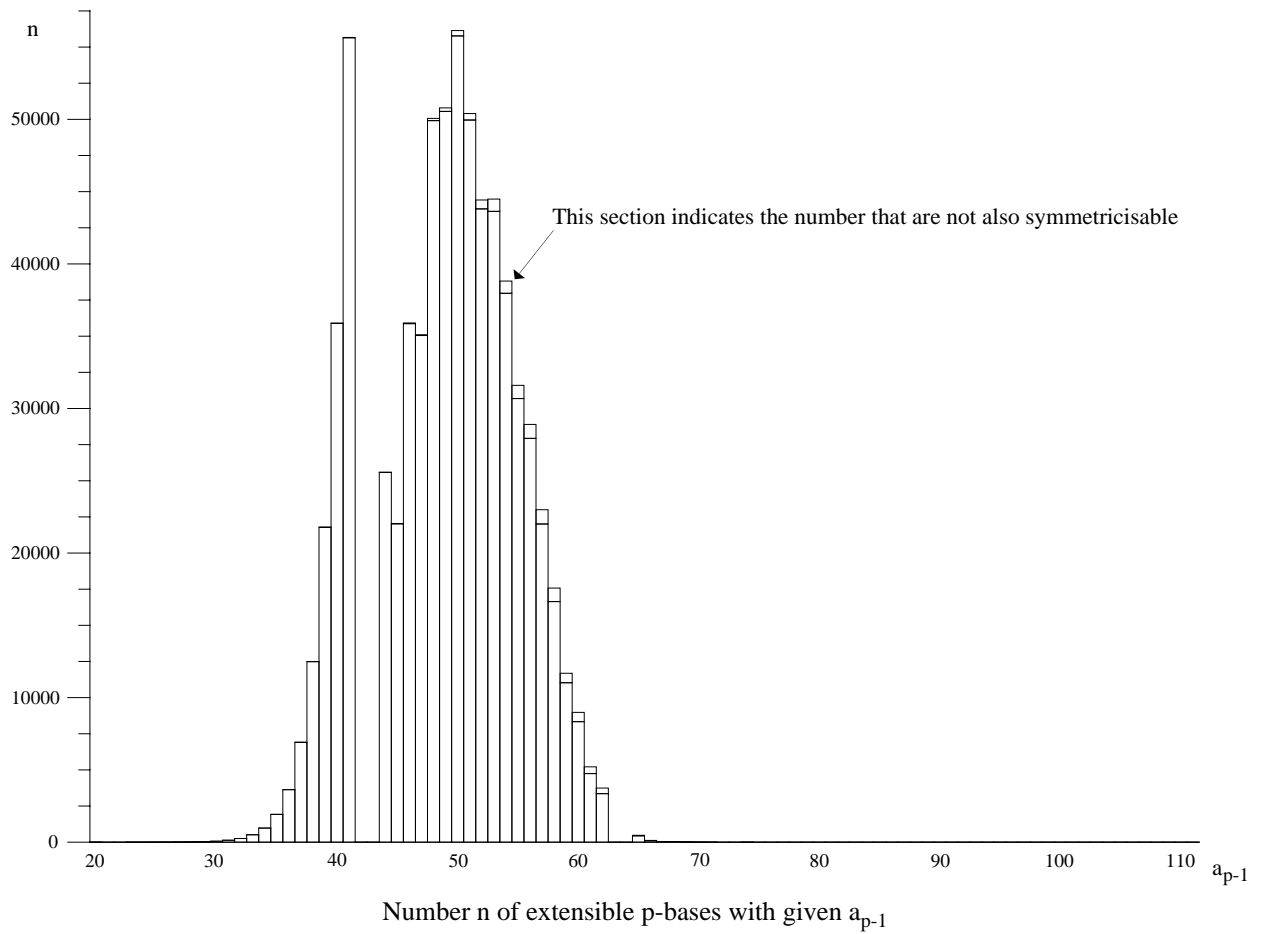

Number n of extensible p-bases with given $a_{p-1}$

Table 13

In Table 14, the plot shows how the chance that a p-basis is also extensible decreases dramatically as $a_{p-1}$ increases.

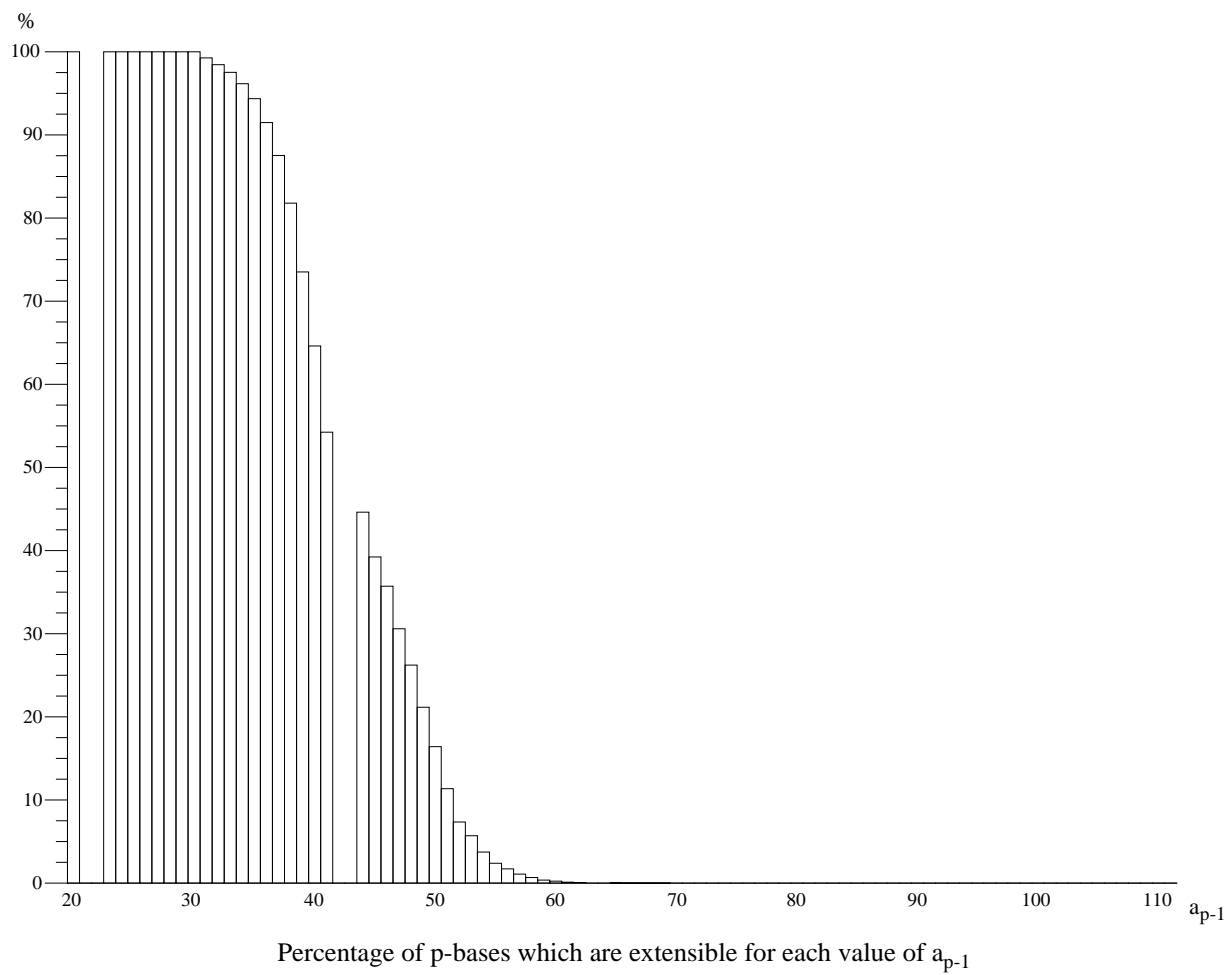

Percentage of p-bases which are extensible for each value of $a_{p-1}$

Table 14

The plot in Table 15 confirms that almost all extensible p-bases are also symmetricisable, even for large values of $a_{p-1}$.

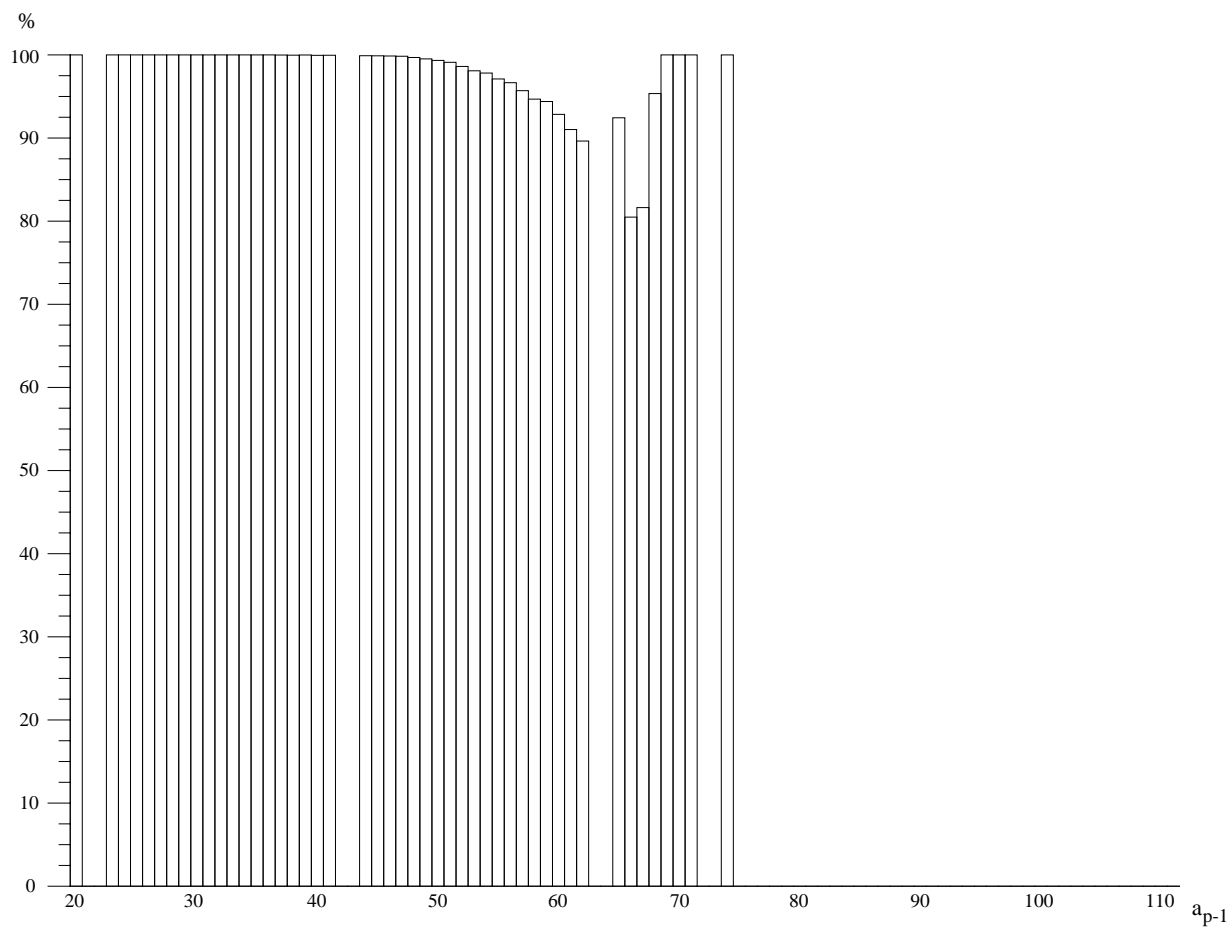

Percentage of extensible p-bases which are also symmetricisable for each value of $a_{p-1}$

Table 15

## 9  Enhanced extensions

We have chosen to extend p-bases in the obvious way (by a sequence of differences p, p, ...) - but are there alternative periodic extensions which have a higher average increment?

A search over the Stohr sequences defined by all p-bases for p<=22 yielded the following three cases where the average increment is as good:

p=18:
$A_{17}$ = {1, 3, 4, 5, 8, 12, 13, 15, 16, 17, 20, 24, 25, 27, 28, 29, 32}  has as Stohr sequence:

```
          62   68   98   104 ...
      30   6    30   6  ...     ie an average of 18
```

p=21:
$A_{20}$ = {1, 3, 4, 5, 6, 9, 14, 15, 17, 18, 19, 20, 23, 28, 29, 31, 32, 33, 34, 37} has as Stohr sequence:

```
          72   79   114  121 ...
      35   7    35   7  ...      ie an average of 21
```

p=22:
$A_{21}$ = {1, 3, 4, 7, 8, 9, 16, 17, 21, 24, 33, 34, 36, 37, 40, 41, 42, 49, 50, 54, 57} has Stohr sequence:

```
          101  112  123  167  178  189 ...
      44   11   11   44   11   11 ...      ie an average of 22
```

No cases where the average increment exceeds p have been found, but I have not been able to prove that this is impossible; on the contrary, the following argument suggests that such bases could exist:

> Suppose the p-basis $A_{p-1}$ is extended by elements $b_0 = a_{p-1}, b_1, b_2, ...$ where the extension sequence is cyclic of length n:
>
> ie    $b_k = b_{k-n} + K$     for some constant K
>
> Choose a sufficiently large value of k (ie one such that $b_k >= 2b_0$), and consider how to generate the K values $b_k <= x < b_{k+n}$.
>
> Because k is sufficiently large, every such generation must include a value $b_i$.
>
> Now suppose:
>
> $b_{k+j} = v_j$ (mod K)      for j = 0 ... n-1
>
> Then the possible values (mod K) that can be generated are:
>
> 1)   $v_0+0, v_0+1, ..., v_0+(p-1)$    - using $b_{k-ln} + a_m$  
>      $v_1+0, v_1+1, ..., v_1+(p-1)$    - using $b_{k+1-ln} + a_m$  
>        ...  
>      $v_{n-1}+0, v_{n-1}+1, ..., v_{n-1}+(p-1)$    - using $b_{k+n-1-ln} + a_m$    for some l, m

2)   $v_i + v_j$                          - using $b_{k+i-ln} + b_{k+j-mn}$     for some l, m

Thus the maximum possible number of different values (mod K) that can be generated is:

   np           - from (1)
   n(n+1)/2     - from (2)

This puts an upper limit on the value of K:

   $K <= np + n(n+1)/2$

Now K=np corresponds to the p, p, ... extension; so this result suggests that complex periodic extensions might be able to improve on this.

## 10  Some further ideas

1)   Is it possible to define a class of "parametric" p-bases which can be used to derive an improved lower bound for n(2,k)?

2)   The extremal basis for k=10 can be thought of as an optimal symmetric sequence whose "tail" has been modified to increase the range; can this idea be applied to other $S^*(p)_k$?